\documentclass[12pt]{amsart}
\usepackage{amsfonts, amssymb, amscd, amsmath}

\newtheorem{thm}{Theorem}[section]
\newtheorem{cor}[thm]{Corollary}
\newtheorem{defn}[thm]{Definition}
\newtheorem{lem}[thm]{Lemma}
\newtheorem{prop}[thm]{Proposition}
\newtheorem{rem}[thm]{Remark}
\numberwithin{equation}{section}
\newtheorem{exl}[equation]{Example}
\def\proof{\textsc{Proof:\ }}
\def\endproof{$\Box$\medskip}

\begin{document}
\title{ Construction of Right Gyrogroup structures on a non-abelian group}
\author{Akhilesh \ Chandra\ Yadav}
\address{Department of Mathematics, M G Kashi Vidyapith, varanasi (INDIA)\\%
\textnormal{akhileshyadav538@gmail.com}}
\today
\begin{abstract}
In this article we gave the notion of particular type of functions on a group termed as class assigned functions . Using  class assigned functions, we have constructed several right gyrogroup structures on a given non-abelian group. \end{abstract}
\maketitle


\section{Introduction}

A gyrogroup is a grouplike structure that is defined in \cite{fog:ung} along with a weaker structure called a left gyrogroup. In \cite{fog:ung1}, Foguel and Ungar   proved that any given group can be turned into a left gyrogroup which is, in turn, a gyrogroup if and only if the given group is central by a $2-Engel$  group. The importance of left gyrogroups stems from the facts that (i) any gyrotransversal groupoid is a left gyrogroup, and (ii) any left gyrogroup is a twisted subgroup in a specifed group. Let $K$  be a group and let $ G$ be the inner product of the group $K$ by $Inn(K)$, where $Inn(K)$ is the inner automorphism group of $K$ whose generic element $\alpha_k$ denotes conjugation by $k\in K$. Then, the diagonal transversal $D$ generated by $K$ (in $G$) is a subset of $G$ given by

$$D = \{(k;\alpha_k)\ \ |\ k \in K\}.$$

Then the  diagonal transversal $D$ is a gyrotransversal  (Theorem 3.2, \cite{fog:ung1}). Indeed this operation induces a binary operation $\odot$ on $K$ so that  the groupoid $(K, \odot)$ is a left gyrogroup (see, \cite{fog:ung1}). The operation $\odot$ on $K$ is given by $x\odot y = x^2 y x^{-1}$. It is convenient for us to study right transversals and right conjugation action. Thus, in our terminology, $(K,\odot)$ is a right gyrogroup, where $x\odot y = y^{-1}xy^2$. In this article, we construct several right gyrogroup structures on a non-abelian group.
\section{Preliminaries}

Let $S$ be a right transversal to the subgroup $H$  in a group $G$. Let $g: S\rightarrow H$ be a map with $g(e) = 1$, where $e$ denotes the identity of $S$ with respect to transversal operation  and $1$ denotes the identity of group $H$.  Consider the right transversal $S_g= \{g(x) x\, |\, x\in S\}$ to $H$. Then  

 \begin{prop}\label{t41}\cite{ylal}
The right quasigroup $(S_g , o')$  is isomorphic to the right quasigroup $(S, o_g)$ where $o_g$ is the operation  given by $xo_g y = x\theta g(y) o y$.
 \end{prop}
The operation $o_g$ is known as the deformed operation on $S$.

\begin{defn}\cite{ylal}\label{gyro}
Let $(S, o)$ be a magma together with a right identity $e$  and together with a right inverse $a'$ for every element $a\in S$ with respect to $e$ in the sense that $ao a'= e$. Then $(S, o)$ will be called a right gyrogroup if \begin{enumerate}
\item  for any $x, y, z\in S$ there exists a unique element 
$x\theta f(y, z)\in S$ such that 
\begin{eqnarray}\label{rq1}
(xoy)oz = x\theta f(y, z)o (yoz),
\end{eqnarray}
\item  the map $f(y,z):S\rightarrow S$ given by $f(y,z) (x)= x\theta f(y, z)$ is an automorphism of $(S, o)$ and 
\item for all $y\in S$,   \begin{eqnarray}\label{rgf}
f(y,  y') = I_S.
\end{eqnarray}
\end{enumerate}
\end{defn}

\begin{defn}\cite{ylal}
  A right transversal $S$ to a subgroup $H$ with identity $e$ in a group  $G$ is called a gyrotransversal if $S = S^{-1}$ and $h^{-1}xh\in S,~\forall~ x\in S $ and $h\in H$. 
\end{defn}  

 \begin{prop}[\cite{ylal}, Representation theorem for right gyrogroups]\label{t52}\
 A right loop $(S, o)$  is a right gyrogroup if and only if it is a  gyrotransversal to $G_S$ in its group extension $G_S S$. 
 \end{prop}

 In \cite{ylal}, we have shown the following:

\begin{lem}\label{t53}\cite{ylal}
 Let $S$ be a gyrotransversal to a subgroup $H$ in a group $G$ and $g$ be a map from $S$ to $H$ such that $g(e) = e$. Then the transversal $S_g = \{g(x) x\ |\ x\in S\}$ is a gyrotransversal if and only if 
 \begin{eqnarray}\label{52}
 g(x^{-1}) &= &g(x)^{-1}
 \end{eqnarray}
and $g$ is equivariant in the sense that 
 \begin{eqnarray}\label{51}
  g(h^{-1} x h) & = & h^{-1}g(x)h 
  \end{eqnarray} 
  for all $x\in S$ and $h\in H.$ 
 \end{lem}

\section{Construction of Right Gyrogroups}
\begin{defn}
Let $G$ be a given non-abelian group.  Define a binary operation $\cdot$ on $G\times G$ as 
$$(a, x)\cdot (b, y) = (ab[b^{-1}xb, y^{-1}], b^{-1}xb o y),$$
where $xoy = y^{-1}xy^2$.
Then $(G\times G, \cdot)$ is a group. The identity of the group $G\times G$
 is $(1, 1)$, where $1$ denotes the identity of group $G$ and the inverse of $(a, x)$ is $(a^{-1}, a^{-1}x^{-1}a)$, for each $(a, x)$ in $G\times G$.
\end{defn}  
One may easily observe that the group $G\times G$ contains $G\times \{1\} (  \simeq G)$ as a subgroup and $\{1\}\times G$ as a gyrotransversal to the subgroup $G\times \{1\}$  in the group $G\times G$(\cite{lal, ylal}).

\begin{prop}
Let $n\in \mathbb{N}$ and $g$ be a map from a group $G$ into itself with $g(1) = 1$ and $g(x )= x^n$ for each $x\in G$. Then $S_g = \{(g(x), x)\, |\, x\in G\} $ is a gyrotransversal to the subgroup $H( = G\times\{1\})$ in $G\times G$.  
\end{prop}

\proof Since $(x^{-1})^n = (x^n)^{-1}$ and $(h^{-1}xh)^n = h^{-1}x^n h$ for all $x, h\in G$. Thus the result follows from Lemma \ref{t53}.
\endproof

Let $G$ be a group. Let  $R$ be a  relation on  $G$ given by $aRb$ if either $a= b$ or $a= b^{-1}$ or $a\sim b$ or $a\sim b^{-1}$, where $\sim$ is a conjugacy relation on group $G$. Then $R$ is an equivalence relation on $G$. The equivalence class of $w\in G$ is denoted by $[w]_R$.

\begin{defn}
Let $G$ be a group. A function $k: G\rightarrow \mathbb{N}\cup\{0\}$ is called a class assigned function on $G$ if $k([1]_R) = 0$ and $k$ is constant on each class $[w]_R$.
\end{defn}

\begin{prop}
Let $k$ be a class assigned function on a free group $F$. Then it determines a function $g: F\rightarrow F$ given by $g(w) = w^{k(w)}$ so that the corresponding transversal 
$S_g= \{(g(w), w)\, |\, w\in F\} $ to the subgroup $H= F\times \{1\}$ in $F\times F$ is a gyrotransversal.
\end{prop}
\proof
By definition of $k$, $k(w ) = k(w^{-1}) = k(u^{-1}wu)$, for all $u, w$ in $F$.
Therefore  $g(w^{-1}) = g(w)^{-1}$ and $g(u^{-1}wu) = (u^{-1}wu)^{k(u^{-1}wu)}= u^{-1}w^{k(u^{-1}wu)}u= u^{-1}w^{k(w)}u= u^{-1}g(w)u.$ 
By Lemma \ref{t53}, $S_g$ is a gyrotransversal.
\endproof

\begin{prop}\label{clsi}
Let $g:F\rightarrow F$ be a map satisfying the property that \\
(i). \,  $ g(1) = 1$,\\
(ii). \,  $g(w^{-1})= g(w)^{-1}$,\\
(iii).\,  $g(u^{-1}wu) = u^{-1}g(w) u$, for all $u, w\in F$. \\
Then $g $ determines a class assigned function $k$  on free group $F$ such that $g(w) = w^{k(w)}$.
\end{prop}

\proof
For $w\in F\setminus \{1\}$, (iii) gives $g(w)= w^{-1} g(w)w$. Thus $g(w)$ belongs to the centralizer of $w$. But $F$ is a free group, therefore $g(w) = w^{k(w)}$, for unique $k(w)\in \mathbb{N}\cup \{0\}$. Then (ii) and (iii) implies that $k(w^{-1}) = k(w)$ and $k(u^{-1}wu) = k(w)$ for all $u,w\in F$. Define $k(1)= 0$, thus $k(w)$ determines a class assigned function $k$.
\endproof

\begin{cor}
There is a bijective correspondence between the set of all class assigned functions and the set of all maps $g$ described as in the Prop. (\ref{clsi}). In particular, there is a bijective correspondence between the set of all class assigned functions and the set of all gyrotransversals to the subgroup $H$ in $F\times F$.
\end{cor}

\begin{thm}
Let $G$ be a non-abelian group and $k$ be a class assigned function on $G$. Then it determines a function $g:G\rightarrow G$ given by $g(w)= w^{k(w)}$ such that $S_g$ is a gyrotransversal to $G\times \{1\}$ in $G\times G$. Define a binary operation $o_k$ on $G$ by 
\begin{eqnarray}\label{d}
x o_k y= y^{-k(y)} x y^{k(y)+1}.
\end{eqnarray}
Then $(G, o_k)$ is a right gyrogroup.  
\end{thm}

\proof
By the property of class assigned function $k$, $g$ satisfies $g(x^{-1})= g(x)^{-1}$ and $g(h^{-1}xh)= h^{-1}g(x)h$ for all $x, h$ in $G$.  Thus $S_g$ is a gyrotransversal to $G\times \{e\}$ in the group $G\times G$ and so deformed operation $o_k$ on $G$ is given by Equation~\ref{d}. Observed that $y^{-1}o_k y= e= y o_k y^{-1}$. By a little computation the map $f(y, z) (x)$ determined by the equation
 $(xoy)oz = f(y, z) (x) o (yoz)$ 
is  $f(y, z) (x)= x \theta \left[y^{k(y)} (yz)^{-k(yo_k z)} z^{k(z)}\right]$,  where $x\theta a= a^{-1}xa$ for all $x,a\in G$. Clearly $f(y, z)$ is an inner automorphism of $(G, o_k)$.
 Also  $f(y,y^{-1})(x)= x=f(y^{-1}, y) (x)$ and so $f(y,y^{-1})=I_G$. Thus, $(G, o_k)$ is a right gyrogroup. \endproof

\begin{rem}
We do not know if its converse true.
\end{rem}

\begin{exl}
Take $G=S_3$, the symmetric group of degree $3$. Let  $k$ be a class assigned function on $S_3$. 

If $k(x)=0$ for all $x$ in $S_3$. Then the corresponding right gyroroup $(S_3, o_k)$ is a group $S_3$ itself.
 
 If  $k(I)=0, k(x)=1$ for all $x\in S_3\setminus\{I\}$. Then the group generated by automorphisms $f(y, z)$ is isomorphic to $[S_3, S_3]=A_3$.

If $k(I)=0$, k((12))=k((2 3))= k((1 3))=1, k((1 2 3))= k((1 3 2))= 2. Then $f\left((12), (132)\right)(x)= x\theta (13)$, $f\left((12), (123)\right)(x)= x\theta (23)$  and $f\left((23), (132)\right)(x)= x\theta (12)$.  Thus,  $(S_3, o_k)$ is a right gyrogroup for which the group generated by automorphisms $f(y, z)$ is  isomorphic to $S_3$.

 The above right gyrogroups are non-isomorphic. Note that  the right gyrogroups determinned by different $k$ may be isomorphic. For, let  $k_1$ and $k_2$ be class assigned maps given by $k_1(x)=0$ for all $x\in S_3$  and $k_2(I)=0$, $k_2(x)=6$ for all $x\in S_3\setminus\{I\}$ respectively. Clearly they are distinct but the corresponding right gyrogroups are the group $S_3$.
 \end{exl}

\begin{cor}
For $n\ge 3$, there are
\begin{eqnarray*}
 \ & \prod_{\sum_{i= 1}^{k} n_i \le n}^{n_i >1}  [n_1, \ldots, n_k] & 
\end{eqnarray*}

 gyrotransversals to the subgroup $S_n$ in the generalized product of group $S_n$. Here $[n_1, \ldots, n_k]$ denotes the least common multiple of $n_1, \ldots, n_k$.
\end{cor}
\proof
Let $p, q\in S_n$. Since $p$ and $p^{-1}$ both have same  cycle structure and $p, q$  are conjugates to each other if and only if they have same cycle structures. Therefore for each $w\in S_n$, the equivalence class $[w]_R$ is same as   the conjugacy class of $w$.  Since  the class assigned function  $k(w)$ has $\left|w\right|$ choices on each class $[w]_R$ and the choices are independent on these classes. Therefore the number of class assigned functions  $k$ will be $\prod \left|w\right|$, where the product runs over each classes with $w\ne 0$. Let $n_1 , n_2, \ldots, n_k$ be natural numbers greater than $1$ such that $n_1 + n_2 + \ldots +n_k\le n$. Then order of permutation with cycle structure $(n_1, n_2, \ldots, n_k)$ is $[n_1, n_2, \ldots, n_k]$, the least common multiple of $n_1, n_2, \ldots, n_k$. Thus,  by above proposition the number of gyrotransversals to $S_n$ in the generalized product $S_n\times S_n$ will be       
$$\prod_{(\sum_{i= 1}^{k} n_i) \le n}^{n_i >1}  [n_1, \ldots, n_k]$$
\endproof

\end{document}